\documentclass[11pt,a4paper]{amsart}
\usepackage{amssymb, amscd, amsmath, enumerate,verbatim}
\usepackage{pb-diagram, amscd}

\input xy
\xyoption{all}

\numberwithin{equation}{subsection}
\newtheorem{theorem}{Theorem}[subsection]
\newtheorem{lemma}[theorem]{Lemma}
\newtheorem{proposition}[theorem]{Proposition}
\newtheorem{corollary}[theorem]{Corollary}

\theoremstyle{definition}

\newtheorem{definition}[theorem]{Definition}
\newtheorem{remark}[theorem]{Remark}

\newcommand{\D}{\mathcal{D}}
\newcommand{\A}{\mathcal{A}}

\newcommand{\K}{\mathcal{K}}

\newcommand{\kd}{\mathcal{KD}}
\newcommand{\kh}{\mathcal{KH}}

\newcommand{\Reg}{\mathbf{Sm}}
\newcommand{\Sch}{\mathbf{Sch}}

\newcommand{\PreSpec}{\mathbf{PreSp}}
\newcommand{\Spec}{\mathbf{Sp}}

\newcommand{\tow}{\mathbf{tow}}

\newcommand{\holim}{\mathrm{holim}}

\newcommand {\lra}{\longrightarrow}

\title[Algebraic $K$-theory and cubical descent]
{Algebraic $K$-theory and cubical descent}

\begin{document}
\begin{large}

\author[Pere Pascual]{Pere Pascual}
\author[Lloren\c{c} Rubi{\'o} Pons]{Lloren\c{c} Rubi{\'o} Pons}
\address
{ Departament de Matem\`{a}tica Aplicada 1\\ Universitat
Polit\`{e}cnica de Catalunya\\Diagonal 647, 08028 Barcelona
(Spain). }
\email{pere.pascual@upc.edu\\llorenc.rubio@upc.edu}

\footnotetext[1]{Partially supported by projects DGCYT
BFM2003-06001 and MT M2006-14575}

\date{\today}

\maketitle

\begin{abstract}
In this note we apply Guill{\'e}n-Navarro descent theorem,
\cite{GN02}, to define a descent variant of the algebraic
$K$-theory of varieties over a field of characteristic zero,
$\mathcal{KD}(X)$, which coincides with $\mathcal{K}(X)$ for
smooth varieties. After a result of Haesemeyer, this new theory is
equivalent to the homotopy algebraic $K$-theory introduced by
Weibel. We also prove that there is a natural weight filtration on
the groups $KH_\ast(X)$.
\end{abstract}

\section{Introduction}
F. Guill{\'e}n and V. Navarro have proved in \cite{GN02} a general
theorem which, in presence of resolution of singularities, permits
to extend some contravariant functors defined on the category of
smooth schemes to the category of all schemes. In this paper we
apply this result to algebraic $K$-theory. More specifically, we
consider the algebraic $K$-theory functor which to a smooth
algebraic variety over a field of characteristic zero $X$
associates the spectrum of the cofibration category of perfect
complexes, $\K(X)$. We apply Guill{\'e}n-Navarro extension criterion
to prove that this functor admits an (essentially unique)
extension to all algebraic varieties, $\kd(X)$, which satisfies a
descent property.

Is is well known that algebraic $K$-theory of schemes does not
satisfies descent. C. Haesmeyer has proved in \cite{H} that the
homotopy algebraic $K$-theory $\kh$, introduced by Weibel in
\cite{W1}, satisfies descent for varieties over a field of
characteristic zero. From the uniqueness of our extension $\kd$
and Haesemeyer's result it follows that, for any variety $X$ over
a field of characteristic zero, the spectra $\kd(X)$ and $\kh(X)$
are weakly equivalent.

Following \cite{GN02} we find also an extension of $\K$ to a
functor with compact support, $\K^c$, which once again by
uniqueness is weakly equivalent to the algebraic $K$-theory with
compact support introduced by Gillet and Soul{\'e} in \cite{GS}.

Moreover, by using the extension theorem in analogy of Guill{\'e}n and
Navarro's paper \cite{GN03}, we are able to prove the existence of
some natural filtrations on the $KD$-groups associated to an
algebraic variety. In fact, the $\kd$-theory of an algebraic
variety $X$ is defined by cubical descent and therefore, if
$X_\bullet$ is a cubical hyperresolution of $X$ (see \cite{GNPP}),
there is a convergent spectral sequence, see proposition
\ref{ssequence},
$$
E_1^{pq}=\bigoplus_{|\alpha|=p+1} K_{q}(X_\alpha)\Rightarrow
KD_{q-p}(X),
$$
where we have written $KD_\ast(X)=\pi_\ast (\kd)$. We prove that
the associated filtration on $KD_{\ast}(X)$ is independent of the
chosen hyperresolution $X_\bullet$ of $X$. In the analogous
situation for compactly supported algebraic $K$-theory we recover
the weight filtration introduced in \cite{GS}. We observe that
Corti{\~n}as, Haesemeyer and Weibel have analyzed in \cite{CHW} the
fiber of the morphism $\K\longrightarrow \kh$ in terms of the
negative cyclic homology functor.

{\em Acknowledgements:} We thank F. Guill{\'e}n, V. Navarro and A.
Roig for many helpful discussions.

\section{The descent theorem of Guill{\'e}n-Navarro}

In this section we recall the main extension theorem proved by
Guill{\'e}n and Navarro and present some corollaries of its proof not
explicitely stated in \cite{GN02}. We also fix some notations.

\subsection{Descent categories} The descent theorem in \cite{GN02} is
stated for functors from the category of smooth varieties to a
cohomological descent category. This kind of category is a
(higher) variation of the classical triangulated categories. We
recall the main features of descent categories and refer to
\cite{GN02}, (1.5.3) and (1.7), for the precise definitions (see
also the proof of proposition \ref{Spec=descent}).

\subsubsection{} For any finite set $S$, the associated {\em cubical} set
$\Box_S$ is the ordered set of non-empty subsets of $S$ and the
{\em augmented cubical} set $\Box^+_S$ is the ordered set of
subsets of $S$, including the empty set. When $S=\{0,1,\cdots ,
n\}$, we simply write $\Box_n$ (respectively, $\Box_n^+$), which
may be identified with the ordered set of $n+1$-tuples $(i_0,\dots
,i_n)$, where $i_k\in\{0,1\}$ such that there is a $k$ with
$i_k\not=0$, and including the $(0,\dots ,0)$ tuple in the
augmented case. We will write $|\alpha | = \sum_{0}^{n}i_k$.

As usual, we will denote by the same symbol the associated
category. Following \cite{GN02}, we denote by $\Pi$ the category
whose objects are finite products of categories $\Box_S$ and whose
morphisms are the functors associated to injective maps in each
component. The objects of $\Pi$ will be called \emph{cubical index
categories}. $\Pi$ is a symmetric monoidal category.

\subsubsection{} Let $\D$ be a category. Given a cubical index category $\Box$,
a $\Box$-\emph{cubical diagram} of $\D$ is a functor
$X:\Box\longrightarrow \D$. We denote by $CoDiag_\Pi \D$ the
category of cubical diagrams of $\D$ (according to \cite{GN02} we
should call these functors {\em cubical codiagrams}, reserving the
term diagram for the contravariant functors
$X:\Box^{op}\longrightarrow \D$): its objects are pairs
$(X,\Box)$, where $X$ is a $\Box$-cubical diagram; a morphism from
the diagram $(X,\Box )$ to the diagram $(Y,\Box ')$ is a functor
$\delta:\Box '\longrightarrow \Box$ together with a natural
transformation $\delta^\ast X= X\circ \delta \Rightarrow Y$.

\subsubsection{}\label{axioms}  A \emph{descent category} is, essentially, a triple
$(\D,E, \mathbf{s})$ given by a cartesian category $\D$ with
initial object $\ast$, a saturated class of morphisms $E$ of $\D$,
called {\em weak equivalences}, and a functor
$$
\mathbf{s} : CoDiag_\Pi \D\longrightarrow \D,
$$
called the {\em simple functor}, which satisfy the following
properties:

\quad $1$. \emph{Product}: for any object $X$ of $\D$, there is a
natural isomorphism $\mathbf{s}_{\Box_0}(X\times \Box_0)\cong X$
and for any $\Box\in \, Ob\Pi$ and any couple of $\Box$-diagrams
$(X,Y)$, the morphism
$$
s_\Box (X\times Y)\longrightarrow s_\Box X\times s_\Box Y,
$$
is an isomorphism.

\quad $2$. \emph{Factorisation}: Let $\Box, \Box'\in\, Ob\Pi$. For
any $\Box\times\Box'$-diagram $X=(X_{\alpha\beta})$, there is an
isomorphism
$$
\mu: s_{\alpha\beta} X_{\alpha\beta}\longrightarrow s_\alpha
s_\beta X_{\alpha\beta},
$$
natural in $X$.

\quad $3$. \emph{Exactness}: Let $f:X\longrightarrow Y$ be a
morphisms of $\Box$-diagrams, $\Box\in\, Ob\Pi$. If for all
$\alpha\in\Box$ the morphism $f_\alpha:X_\alpha\longrightarrow
Y_\alpha$ is a weak equivalence (i.e. it is in $E$), then the
morphism $s_\Box f : s_\Box X\longrightarrow s_\Box Y$ is a weak
equivalence.

\quad $4$. \emph{Acyclicity criterium}: Let $f:X_1\longrightarrow
X_0$ be a morphism of $\D$. Then, $f$ is a weak equivalence if and
only if the simple of the $\Box_1$-diagram
$$
\ast \longrightarrow X_0 \stackrel{f}{\longleftarrow }X_1
$$
is acyclic, that is, it is weakly equivalent to the final object
of $\D$.

The acyclicity criterium has to be verified also for higher
cubical diagrams, \cite{GN02}. More specifically, let $X^+$ be a
$\Box_n^+$-diagram in $\D$ and denote by $X$ the cubical diagram
obtained from $X^+$ by restriction to $\Box_n$. Then the
acyclicity criterium takes the following form (see property
(CD8)$^{op}$ of definition $(1.5.3)$ of \cite{GN02}):

\quad $4'$. \emph{Acyclicity criterium}: The augmentation morphism
$\lambda_\varepsilon: X_0\longrightarrow s_\Box X$ is a weak
equivalence if and only if the canonical morphism
$\ast\longrightarrow s_{\Box^+}X^+$ is a weak equivalence.

We remark that the transformations $\mu$ and $\lambda$ of
properties $2$ and $4'$ are, in fact, part of the data of a
descent structure.

\subsubsection{} The categories of complexes give the basic examples of descent
categories: if $\A$ is an abelian category, the category of
bounded below cochain complexes ${\mathbf C}^\ast(\A)$, with the
class of quasi-isomorphisms as weak equivalences and the total
functor of a multicomplex as simple functor, is a descent
category. See \cite{GN02} for other examples.

\subsection{Guill{\'e}n-Navarro theorem}
Let $k$ be a field of characteristic zero. We denote by $\Sch(k)$
the category of reduced separated schemes of finite type over $k$,
simply called \emph{algebraic varieties}, and by $\Reg(k)$ the
category of smooth varieties.

\subsubsection{}\label{square}
Let
$$
\begin{CD}
\widetilde{Y} @>j>> \widetilde{X}\\\
@VgVV @VVfV \\
 Y @>i>> X
 \end{CD}
$$
be a cartesian diagram of schemes, which we may consider as a
$\Box^+_1$-diagram. We say that it is an {\em acyclic square} if
$i$ is a closed immersion, $f$ is a proper morphism and the
induced morphism
$\widetilde{X}\setminus\widetilde{Y}\longrightarrow X\setminus Y$
is an isomorphism.

We say that an acyclic square is an {\em elementary acyclic
square} if all schemes in the diagram are irreducible and smooth,
and $f$ is the blow-up of $X$ along $Y$.

\begin{theorem}\label{extensio}{\rm (\cite{GN02}, (2.1.5))}
Let $\D$ be a cohomological descent category and
$$
G: \Reg(k)\longrightarrow {\rm Ho} {\D}
$$
a contravariant $\Phi$-rectified functor satisfying the following
conditions:
\begin{itemize}
\item[(F1)] $G(\emptyset)=0$, and the canonical morphism
$G(X\sqcup Y)\longrightarrow G(X)\times G(Y)$ is an isomorphism,
\item[(F2)] if $X_\bullet$ is an elementary acyclic square in
$\Reg(k)$, then $sG(X_\bullet)$ is acyclic.
\end{itemize}
Then there is an extension of $G$ to a $\Phi$-rectified functor
$$
GD:\Sch(k)\longrightarrow {\rm Ho}{\D}
$$
which satisfies the {\em descent condition}
\begin{itemize}
\item[(D)] if $X_\bullet$ is an acyclic square in $\Sch(k)$,
$sGD(X_\bullet)$ is acyclic.
\end{itemize}
Moreover, this extension is essentially unique: if $G'$ is another
extension of $G$ verifying the descent property (D), then there is
a uniquely determined isomorphism of $\Phi$-rectified functors
$GD\Rightarrow G'$.
\end{theorem}

We will say that the functor $GD$ has been obtained from $G$ by
cubical descent.

The proof of Guill{\'e}n-Navarro's theorem gives more than stated
above. In fact, if $X$ is an algebraic variety and
$X_\bullet\longrightarrow X$ is any cubical hyperresolution, see
\cite{GNPP}, it is proved in \cite{GN02} that, under the
hypothesis of the theorem,
$$
GD(X)= sG(X_\bullet ),
$$
gives a well defined functor from $\Sch(k)$ to $\mathrm{Ho}\D$,
independent of the chosen hyperresolution $X_\bullet$. From this
explicit presentation we deduce easely some more properties of the
descent extension $GD$.

\begin{proposition}\label{GaGD}
Suppose that the functor $G$ in theorem {\rm\ref{extensio}} is
already defined for all varieties, that is, we have $G:\Sch(k)\lra
{\rm Ho}{\D}$, and satisfies {\rm (F1)} and {\rm (F2)}. Then there
is a natural transformation of $\Phi$-rectified functors
$G\Rightarrow GD$.
\end{proposition}
\begin{proof}
Let $X$ be a variety and $X_\bullet$ a cubical hyperresolution of
$X$, indexed by a cubical set $\Box$. Taking the simple of the
morphism of cubical diagrams $G(X\times\Box)\lra G(X_\bullet)$ we
get the morphism $G(X)=sG(X\times\Box)\lra sG(X_\bullet)=GD(X)$.
\end{proof}

Looking at the construction and properties of cubical
hyperresolutions, it may be proved that the extended functor $GD$
inherits many properties of the functor $G$ over the smooth
varieties. As an example, and in view of their interest in
algebraic \emph{K}-theory, let us remark the two properties
inclosed in the following proposition.

\begin{proposition}\label{propietatsGD}
Consider the hypothesis of theorem {\rm\ref{extensio}}.
\begin{itemize}
\item[(1)] Suppose that $G$ is homotopy invariant, i.e. for any
smooth variety $X$ the projection $X\times \mathbb{A}^1\lra X$
induces an isomorphism $G(X)\cong G(X\times\mathbb{A}^1)$. Then
$GD$ is homotopy invariant: for any variety $X$, there is an
isomorphism $GD(X)\cong GD(X\times\mathbb{A}^1)$. \item[(2)]
Suppose that $G$ satisfies the Mayer-Vietoris, i.e. for any smooth
variety $X$ and any open decomposition $X=U\cup V$ the square,
induced by inclusions,
$$
\begin{CD}
G(X) @>>> G(U)\\\
@VVV @VVV \\
G(V)@>>> G(U\cap V)
\end{CD}
$$
is acyclic in $\D$. Then $GD$ satisfies Mayer-Vietoris for all
varieties.
\end{itemize}
\end{proposition}
\begin{proof}
Given $X$ an algebraic variety we fix $X_\bullet $, a cubical
hyperresolution of $X$.

(1) By the definition of $GD$ and the homotopy invariance of $G$
we have a sequence of weak equivalences
$$
GD(X)\cong sG(X_\bullet)\cong sG(X_\bullet\times
\mathbb{A}^1)\cong GD(X\times \mathbb{A}^1),
$$
so the proof follows.

(2) By the definition of cubical hyperresolutions (see
\cite{GNPP}), the restrictions of $X_\bullet$ to $U,V$ and $U\cap
V$ give hyperresolutions of these varieties. Let's denote by
$U_\bullet ,V_\bullet $ and $(U\cap V)_\bullet$, respectively,
these restrictions. By construction, for any index $\alpha$ we
have an open decomposition $X_\alpha = U_\alpha\cup V_\alpha$ with
$U_\alpha \cap V_\alpha = (U\cap V)_\alpha$, so from the
Mayer-Vietoris property for $G$ on the category of smooth schemes,
we deduce that the morphisms
$$
G(X_\alpha ) {\longrightarrow } s\left( \begin{CD} @.
G(U_\alpha)\\
@. @VVV\\
G(V_\alpha ) @>>> G((U\cap V)_\alpha )
\end{CD}\right)
$$
are weak equivalences for any $\alpha$. By the exactness property
of descent categories, we have that
$$
sG(X_\bullet ) {\longrightarrow } s_\alpha s\left(
\begin{CD} @.
G(U_\alpha)\\
@. @VVV\\
G(V_\alpha ) @>>> G((U\cap V)_\alpha )
\end{CD}\right)
$$
is also a weak equivalence. But, by the factorization axiom of
descent categories, the simple on the right is weak equivalent to
$$
s\left( \begin{CD} @.
sG(U_\bullet )\\
@. @VVV\\
sG(V_\bullet ) @>>> sG((U\cap V)_\bullet )
\end{CD}\right)
$$
So, taking into account the definition of $GD$ we finally deduce
that the morphism $GD(X)\lra s(GD(U)\longleftarrow GD(U\cap V)\lra
GD(V))$ is a weak equivalence, hence tha Mayer-Vietoris property
for open sets follows.
\end{proof}

\subsection{Extension with compact support} In \cite{GN02},
the authors present some variations on the main theorem. In
particular, they prove in \cite{GN02}, $(2.2.2)$, that with the
same hypothesis of theorem \ref{extensio} there is an extension
$G^c$ of $G$ with compact support: if $\Sch_c(k)$ denotes the
category of varieties and proper morphisms, there is an extension
of $G$ to a $\Phi$-rectified functor
$$
G^c:\Sch_c(k)\longrightarrow {\rm Ho}{\D}
$$
which satisfies the descent property (D) and, moreover,
\begin{itemize}
\item[($D_c$)] if $Y$ is a subvariety of $X$, then there is a
natural isomorphism
$$
G^c(X-Y) \cong s_{\Box^+_0}(G^c(X)\longrightarrow G^c(Y)).
$$
\end{itemize}

\section{The descent category of Spectra}

In this section we prove that the category of $\Omega$-spectra,
with the homotopy limit as a simple functor, is a (cohomological)
descent category in the sense of \cite{GN02}.

\subsection{Fibrant spectra} We will work in the category of
fibrant spectra of simplicial sets. Our main references will be
the paper by Bousfield-Friedlander \cite{BF} and section $5$ of
Thomason's \cite{T80}.

Recall that a {\em prespectrum} is a sequence of pointed
simplicial sets $X_n$, $n\geq 0$, together with structure maps
$\Sigma X_n\longrightarrow X_{n+1}$, where for a pointed
simplicial set $K$, $\Sigma K=S^1\wedge K$. A prespectrum $X$ is a
{\em fibrant spectrum}, also called $\Omega$-\emph{spectrum}, if
each $X_n$ is a fibrant simplicial set and the maps
$X_n\longrightarrow \Omega X_{n+1}$, obtained by adjunction of the
structure maps, are weak equivalences. Morphisms between
preespectra and between fibrant spectra are defined as maps in
each degree that commute with the structure maps. We denote by
$\PreSpec$ and $\Spec$ the categories of prespectra and fibrant
spectra, respectively.

The homotopy groups of a prespectrum $X$ are defined by the direct
limit
$$
\pi_k(X) = \varinjlim \pi_{k+n}(X_n),\quad k\in\mathbb{Z},
$$
so that if $X$ is a fibrant spectrum, $\pi_k(X) = \pi_{k+n}(X_n)$
for $k+n\geq 0$, and, more specifically, for $k\geq 0$, $\pi_k(X)
= \pi_{k}(X_0)$. A map $f:X\longrightarrow Y$ of prespectra is a
{\em weak equivalence} if it induces an isomorphism on homotopy
groups. In this way, a map of fibrant spectra is a weak
equivalence if and only if it induces weak equivalences in each
degree.

\subsection{Homotopy limit}
Let $X$ be a functor from an index category $I$ to $\Spec$. The
homotopy limit spaces $\holim X_n$, $n\geq 0$, in the sense of
Bousfield-Kan, \cite{BK}, chapter XI, define a fibrant spectrum,
$\holim X$, see \cite{T80}, 5.6. In fact, one can see that
$\PreSpec$ has a structure of simplicial closed model category,
see \cite{S}, so that we can apply the general theory of homotopy
limits for theses categories, \cite{H}.

The main properties we need of homotopy limits between fibrant
spectra are:
\begin{itemize}
 \item[$(i)$] Functoriality and exactness on fibrant spectra: Let
$f:X\lra Y$ be a morphism of $I$-diagrams spectra. Then, there is
a natural morphism $\holim f : \holim X\lra \holim Y$. If for each
$\alpha\in I$ the morphism $f_\alpha:X_\alpha\lra Y_\alpha$ is a
weak equivalence, then $\holim f$ is a weak equivalence.
\item[$(ii)$] Functoriality on the index category and cofinality
theorem: Given a functor $\delta:I\lra J$ and a diagram $X:J\lra
\Spec$, there is a natural map $\holim _J X \lra \holim_I
\delta^*X$, where $\delta^*X =X\circ \delta$. If $\delta$ is left
cofinal, this morphism is a weak equivalence. \item[$(iii)$] For
any diagram $X:I\lra \Spec$, there is a natural map $\lim X\lra
\holim X$.

\end{itemize}

\subsubsection{} For a cubical diagram of spectra $X:\Box\longrightarrow \Spec$
we define the {\em simple spectrum of} $X$ as the homotopy limit
$$
s_\Box (X) = \holim_\Box X.
$$
For a fixed cubical category $\Box$, $s_\Box$ defines a functor
$s_\Box: CoDiag_\Box\Spec \longrightarrow \Spec$, and by the
functoriality of the homotopy limit with respect to the index
category $\Box$, we obtain a functor
$$
{\mathbf s}: CoDiag_\Pi\Spec \longrightarrow \Spec .
$$
\subsubsection{}\label{def+} Following \cite{GN02} (1.4.3), we extend the functor
$s$ to augmented cubical diagrams by using the cone construction.
For instance, if $f:X\longrightarrow Y$ is a $\Box_0^+$-diagram of
spectra, that is to say, a morphism, it follows from {\em loc.
cit.} that
$$
s_{\Box_0^+} (f) = s_{\Box_1} (X
\stackrel{f}{\longrightarrow} Y \longleftarrow \ast ),
$$
which is weakly equivalent to the homotopy fiber of $f$.

Take an isomorphism $\Box_n^+ \cong \Box_0^+ \times \Box_{n-1}^+$.
As the cone construction respects this product structure, we find
$$
s_{\Box_n^+} X^+ = s_{\Box_{n-1}^+}(s_{\Box_0^+}X^+),
$$
that is, by viewing $X^+$ as a morphism of two
$\Box^+_{n-1}$-diagrams, $f:X^+_0\lra X^+_1$, the simple spectrum
associated to $X^+$ is obtained as the simple of the
$\Box_{n-1}^+$-cubical diagram which in each degree $\alpha$ has
the homotopy fiber of $f_\alpha$. As a consequence, the simple
spectrum $s_{\Box_n^+} X^+$ is isomorphic to the total fiber space
of $X^+$ as defined by Goodwillie in \cite{G}.

\subsubsection{} If $X^+$ is a $\Box_n^+$-diagram and $X$ denotes its
restriction to $\Box_n$, it follows from the general properties of
homotopy limits outlined above that there is a natural map
$X_0\lra \holim X$. As a consequence of \cite{G}, 1.1.b (compare
also with \cite{P}, proposition (3.3), for a similar situation),
we obtain:

\begin{proposition}\label{good}
Let $X^+:\Box_n^+\longrightarrow \Spec$ be an augmented cubical
diagram of spectra and $X$ the cubical diagram obtained by
restriction to $\Box_n$. The simple $s_{\Box_n^+} X^+$ is
isomorphic to the homotopy fiber of the morphism
$X_0\longrightarrow s_\Box X=\holim X$.\hfill $\Box$
\end{proposition}

Denote by $\ast$ the initial object of $\Spec$. The following
corollary relates the simple of a cubical diagram with the simple
of an augmented diagram.

\begin{corollary}
Let $X:\Box_n\longrightarrow \Spec$ be a cubical diagram of
spectra and let $\widetilde X$ the augmented cubical diagram
obtained from $X$ by adding $X_0=\ast$. Then,
$$
s_{\Box_n^+}\widetilde X = \Omega s_{\Box_n} X.
$$
\end{corollary}

We also deduce the following result, which will be used later:

\begin{corollary}\label{ssequence1}
Let $X_\bullet$ be a $\Box_n$-diagram of spectra. Then, there is a
convergent spectral sequence
$$
E_1^{pq} = \bigoplus_{|\alpha |=p+1}\pi_q(X_\alpha)\Longrightarrow
\pi_{q-p}(s_{\Box_n}X_\bullet).
$$
\end{corollary}
\begin{proof}
Consider the cubical diagrams $F^pX_\bullet$ defined by
$$
(F^pX_\bullet)_\alpha = \left\{
\begin{array}{cc}
  X_\alpha , & \mbox{if}\quad |\alpha| \leq p+1,\\
  \ast , & \mbox{if}\quad |\alpha| > p+1. \\
\end{array}
\right.
$$
Observe that $F^{-1}X_\bullet$ is the constant diagram defined by
$\ast$ and that $F^{n}X_\bullet = X_\bullet$. We obtain a sequence
of cubical diagrams
$$
F^{n}X_\bullet\lra F^{n-1}X_\bullet\lra \dots \lra
F^0X_\bullet\lra \ast
$$
which is a degreewise sequence of fibrations of spectra. Hence,
taking homotopy limits there is a sequence of fibrations
$$
s_\Box(F^{n}X_\bullet )\lra s_\Box(F^{n-1}X_\bullet )\lra \dots
\lra s_\Box(F^0X_\bullet )\lra \ast .
$$
The Bousfield-Kan spectral sequence associated to the tower of
fibrations obtained by adjoining identities from the left
converges to the homotopy of $s_\Box (F^{n} X_\bullet )=s_\Box
X_\bullet$. The $E_1$ terms are $E_1^{pq}=\pi_{q-p}(sGr^p
X_\bullet)$, where $Gr^pX_\bullet$ is the $\Box$-diagram obtained
degree wise as the fibers of the morphism $s_\Box(F^pX_\bullet
)\lra s_\Box(F^{p-1}X_\bullet )$. But, reasoning as in the proof
of proposition $(3.3)$ of \cite{P}, for these diagrams we have
$$
s_{\Box_n} Gr^pX_\bullet = \prod_{|\alpha|=p+1}\Omega^pX_\alpha ,
$$
hence it follows that
$$
E_1^{pq} = \pi_{q-p}(\prod_{|\alpha|=p+1}\Omega^pX_\alpha) =
\bigoplus_{|\alpha |=p+1}\pi_q(X_\alpha).
$$
Convergence is a consequence of lemma 5.48 of \cite{T80}.
\end{proof}

\subsubsection{}
We say that an augmented cubical diagram of spectra $X^+$ is
\emph{acyclic} if the canonical morphism $\ast\longrightarrow
s_{\Box^+}X^+$ is a weak equivalence. The acyclic diagrams are
also called \emph{homotopy cartesian} diagrams, see \cite{G} and
\cite{W1}. From proposition \ref{good} it follows immediately (see
also \cite{W1}, proposition 1.1):

\begin{corollary}\label{cd8}
Let $X^+:\Box_n^+\longrightarrow \Spec$ be an augmented cubical
diagram of spectra and $X$ the cubical diagram obtained by
restriction to $\Box_n$. Then $X^+$ is acyclic if and only if the
natural morphism $X_0\longrightarrow \holim X$ is a weak
equivalence. \hfill $\Box$
\end{corollary}

\begin{remark}
Observe that for $n=2$ this result reduces to the well known fact
that a square of fibrant spectra
$$
\begin{CD}
X @>>> Y\\
@VVV @VVhV\\
X' @>k>> Y'
\end{CD}
$$
is acyclic (or homotopy cartesian) if and only if the natural map
from $X$ to the homotopy limit of $X'\stackrel{k}{\longrightarrow}
Y'\stackrel{h}{\longleftarrow} Y$ is a weak equivalence.
\end{remark}

\subsubsection{}
After the remarks above, on $\Spec$ we have a class of weak
equivalences and a simple functor $ {\mathbf s}: CoDiag_\Pi\Spec
\longrightarrow \Spec$. According to \cite{GN02}, d{\'e}finition
($1.7.1$), to have a (cohomological) descent category on $\Spec$
we also need the following data:
\begin{itemize}
\item[$(i)$] a natural transformation $\mu : s_\Box\circ s_{\Box
'}\Rightarrow s_{\Box\times\Box '}$, \item[$(ii)$] a natural
transformation $\lambda_\Box : id_{\Spec}\Rightarrow s_\Box\circ
i_\Box$,
\end{itemize}
in such a way that $(s,\mu ,\lambda_0):\Pi\lra CoReal_\Pi\Spec$
defines a comonoidal quasi-strict functor, see \emph{loc.cit.} As
the homotopy limit is the end of a functor, by Fubini theorem (see
\cite{M})  there is a natural transformation
$$
\mu : s_\Box\circ s_{\Box '}\Rightarrow s_{\Box\times\Box '},
$$
such that for any diagram $X$,  $\mu_X$ is an isomorphism. As for
$\lambda$, recall that $s_\Box(X\times \Box)$ is the function
space from the classifying space of the index category $\Box $ to
$X$, so one defines
$$
\lambda_\Box (X) : X\lra s_\Box(X\times \Box),
$$
by constant functions.

\begin{proposition}\label{Spec=descent}
The category of fibrant spectra $\Spec$ with weak homotopy
equivalences as weak equivalences and the homotopy limit $\holim$
as simple functor for cubic diagrams, and the natural
transformations $\mu ,\lambda$ defined above, is a cohomological
descent category.
\end{proposition}
\begin{proof}
The actual definition of cohomological descent category consist of
8 axioms, which are dual to the axioms (CD1)-(CD8) of \cite{GN02},
definition $(1.5.3)$, (see also their $(1.7)$). Much of them are
immediate from the definitions and the properties of homotopy
limits, so we comment the four axioms summarized in section
\ref{axioms}, (see also \cite{R} for an extension of this result
to stable simplicial model categories).

It is clear from the definitions that $\Spec$ is a cartesian
category with initial object $\ast$.

$1$. \emph{Product}: since the homotopy limit is an end, it is
compatible with products, so for any $\Box$-diagrams $X,Y$ of
$\Spec$ there is a natural isomorphism
$$
s_\Box(X\times Y) \cong s_\Box(X) \times s_\Box(Y).
$$

$2$. \emph{Factorisation}: also because of the Fubini theorem for
ends, if $X$ is a $\Box\times\Box'$-diagram, there are natural
isomorphism
$$
s_\Box s_{\Box'} X_{\alpha\beta}\cong
s_{\Box\times\Box'}X_{\alpha\beta}\cong s_{\Box '}s_\Box
X_{\alpha\beta},
$$
see \cite{T80}, lemma 5.7.

$3$. \emph{Exactness}: If $f:X\longrightarrow Y$ is a morphism of
$\Box$-diagrams in $\Spec$ such that for any $\alpha\in\Box$ the
morphism $f_\alpha$ is a weak equivalence, then $s_\Box f:s_\Box
X\longrightarrow s_\Box Y$is a weak equivalence, since the homotpy
limit preserves weak equivalences between fibrant spectra, see
\cite{T80}, 5.5. Observe that this property is not true for
prespectra.

$4'$ \emph{Acyclicity criterium}: this is exactly the result of
corollary \ref{cd8}.
\end{proof}

\section{Descent algebraic $K$-theory}

\subsection{} Let $X$ be a noetherian separated scheme, we denote
by $\K(X)$ the $K$-spectrum associated to the category of perfect
complexes on $X$, see \cite{TT}, definition 3.1. It defines a
contravariant functor from the category of noetherian separated
schemes to the category of spectra $\Spec$ (\cite{TT}, 3.14).
Moreover, it is a covariant functor for perfect projective maps
and for proper flat morphisms, (\cite{TT}, 3.16).

\begin{theorem}\label{teorema1}
Let $k$ be a field of characteristic zero. The functor
$$
\K : \Reg(k)\longrightarrow {\rm Ho}\Spec ,
$$
admits a unique extension, up to unique isomorphism of
$\Phi$-rectified functors, to a functor
$$
\kd : \Sch(k)\longrightarrow {\rm Ho} \Spec
$$
such that satisfies the descent property (D):
\begin{itemize}
\item[(D)] if $X_\bullet$ is an acyclic square in $\Sch(k)$,
$s\kd(X_\bullet)$ is acyclic.
\end{itemize}
\end{theorem}
\begin{proof}
By proposition \ref{Spec=descent}, we know that $\Spec$ is a
descent category. So, in order to apply  Guill{\'e}n-Navarro descent
theorem \ref{extensio} we have to verify properties $(F1), (F2)$.
The first one is immediate, while $(F2)$ follows from Thomason's
calculation in \cite{T93} of the algebraic $K$-theory of a blow up
along a regularly immersed subscheme, as has been observed by many
authors (see, for example,  in \cite{H},  \cite{GS} and
\cite{CHSW}).

In the context of cubical spectra we propose the following
presentation of property $(F2)$. Consider an elementary acyclic
square as in \ref{square} and the square of spectra obtained by
application of the algebraic $K$ functor
$$
\begin{CD}
\K(X) @>i^\ast >> \K(Y)\\
@Vf^\ast VV @VVg^\ast V\\
\K(\widetilde{X}) @>j^\ast >> \K(\widetilde{Y})
\end{CD}
$$
We have to prove that this square is an acyclic square of spectra.
If $N$ is the conormal bundle of $Y$ in $X$, then $\widetilde
Y={\mathbb P}(N)$, so the morphism
$$
\Psi: \prod^{d}\K(Y)\longrightarrow \K(\widetilde{Y}),
$$
induced by the functor which is defined on a sequence of perfect
complexes by
$$
(E_0,\dots,E_{d-1})\mapsto \bigoplus_{i=0}^{d-1} {\mathcal
{O}}_{{\mathbb P}(N)}(-i)\otimes Lg^* E_i,
$$
is a weak equivalence, see \cite{TT}, theorem 4.1, and also
\cite{T91}.

For the blown up variety $\widetilde X$, it has been proved by
Thomason, see \cite{T93}, th{\'e}or{\`e}me 2.1, that the morphism
$$
\Phi: \K(X)\times \prod^{d-1}\K(Y) \longrightarrow \K(\widetilde
X)
$$
which is induced by the functor on perfect complexes given by
$$
(F,E_1,\dots ,E_{d-1})\mapsto f^*F \oplus\bigoplus_{i=1}^{d-1}
j_\ast({\mathcal {O}}_{{\mathbb P}(N)}(-i)\otimes Lg^* E_i),
$$
is also a weak equivalence.

Define $j': \K(X)\times \prod^{d-1}\K(Y)\longrightarrow
\prod^{d}\K(Y)$ componentwise by $g^*i^*$ on the first component
and the morphism given by multiplication by $\lambda_{-1}(N)$ in
the $Y$-components. After the self-intersection formula,
\cite{T93}, (3.1.4), the diagram
$$
\begin{CD}
 \K(X)\times \prod^{d-1}\K(Y) @>\Phi>>\K(\widetilde X)\\
@Vj' VV @VVj^\ast V\\
\prod^{d}\K(Y)@>\Psi>> \K(\widetilde{Y})
\end{CD}
$$
is commutative. Since $\Phi, \Psi$ are weak equivalences, it is an
acyclic diagram.

Consider now the augmented commutative cubical diagram
$$
\xymatrix{
                   & \K(X) \ar@{=}[ld] \ar[rr]
                   \ar@{-}[d]
                   &            & \K(X)\times \prod^{d-1} \K(Y) \ar[dl]_\Phi \ar[dd]^{j'}\\
\K(X) \ar[dd]^{i^\ast} \ar[rr]^{f^\ast} & \ar[d]
& \K(\widetilde X) \ar[dd]_{j^\ast} &\\
                   & \K(Y)     \ar@{=}[ld] \ar@{-}[r]  &    \ar[r]    &
                   \prod^d \K(Y) \ar[dl]_\Psi \\
\K(Y)\ar[rr]^{g^\ast}         &                              &
\K(\widetilde Y) & }
$$
where the horizontal back arrows are the inclusion on the first
factor.

As the right and left side squares are acyclic, it follows from
the definition in \ref{def+} that it is an acyclic cubical
diagram. But the back square is acyclic because the two horizontal
morphisms have the same cofiber, so the front square must be
acyclic, which is what has to be proved.
\end{proof}

For a $k$-variety $X$, we will denote by $KD_\ast(X)$ the homotopy
groups of $\kd(X)$, $$KD_\ast(X):=\pi_\ast(\kd(X)).$$ The descent
property (D) gives rise to exact sequences:

\begin{corollary}\label{exacta0}
Let $X_\bullet$ an acyclic square in $\Sch(k)$. Then there is an
exact sequence
$$
\dots \lra KD_n(X)\stackrel{f^\ast -i^\ast}{\lra }KD_n(\widetilde
X)\oplus KD_n( Y)\stackrel{j^\ast + g^\ast}{\lra }KD_n({\widetilde
Y})\stackrel{\delta}{\lra }KD_{n-1}(X)\lra \dots
$$
\end{corollary}

More generally, if $X$ is a $k$-variety, then $\kd (X)$ is defined
as the simple of the cubical diagram of spectra $\K (X_\bullet)$,
where $X_\bullet$ is a cubical hyperresolution, so from
proposition \ref{ssequence1} we deduce:
\begin{proposition}\label{ssequence}
Let $k$ be a field of characteristic zero and $X$ be an algebraic
$k$-variety. Let $X_\bullet$ be a cubical hyperresolution of $X$.
Then, there is a convergent spectral sequence
$$
E_1^{pq} = \bigoplus_{|\alpha|=p+1} K_q(X_\alpha) \Longrightarrow
KD_{q-p}(X).
$$
\end{proposition}

If $X$ is of dimension $d$, we can take cubical hyperresolutions
of size $d$ (see \cite{GNPP}, I.2.15), so it follows:

\begin{corollary}
Let $k$ be a field of characteristic zero and $X$ be an algebraic
$k$-variety of dimension $d$. Then,
$$
KD_n(X) = 0, \qquad n< -d.
$$
\end{corollary}

\subsection{Some properties of $\kd$}
As explained in section 1, $\kd$ inherites many properties of the
algebraic $K$-theory of smooth schemes. For example, from
proposition \ref{propietatsGD} and the properties of homotopy
invariance and Mayer-Vietoris for the $K$-theory of smooth
schemes, (see \cite{Q}), we deduce immediately:
\begin{proposition}
The descent $\kd$-theory satisfies:
\begin{itemize}
\item[{\rm (1)}] $\kd$ is homotopy invariant, that is, for any
variety $X$ the projection $X\times \mathbb{A}^1\longrightarrow X$
induces a weak equivalence $ \kd(X)\cong \kd(X\times
\mathbb{A}^1).$ \item[\rm{(2)}] $\kd$ has the Mayer-Vietoris
property, that is, if $X=U\cup V$, with $U,V$ open sets, then the
square
$$
\begin{CD}
{\mathcal KD}(X) @>>>  {\mathcal KD}(U)\\
@VVV @VVV\\
{\mathcal KD}(V) @>>> {\mathcal KD}(U\cap V)
\end{CD}
$$
is homotopy cartesian.\hfill{$\Box$}
\end{itemize}
\end{proposition}

One may prove in a similar maner that $\kd$ satisfies the
fundamental Bass theorem. Also following \cite{T80}, \cite{W1} one
can prove the existence of a Brown-Gersten type spectral sequence:
$$
E_2^{pq} = H^p(X,\widetilde{KD}_{-q})\Rightarrow KD_{-p-q}(X).
$$
where $\widetilde{KD}_{\ast}$ stands for the sheaf in the Zariski
topology associated to the presheaf ${KD}_{\ast}$.

\subsection{Equivalence with homotopy algebraic  $K$-theory}
In \cite{H}, theorem 3.5, Haesemeyer has proved that the homotopy
algebraic $K$-theory $\kh$ of an algebraic variety $X$ defined by
Weibel in \cite{W1} (see also \cite{TT}), satisfies the descent
axiom (D). As the $KH$-theory coincides with $K$-theory for smooth
varieties, we can apply the uniqueness property of the extension
theorem \ref{extensio} to obtain:

\begin{corollary}\label{kh=kd}
Let $X$ be an algebraic variety over a field of characteristic
zero. There is a natural morphism $\kd(X)\lra \kh(X)$, in
$\mathrm{Ho}\Spec$, which is a weak equivalence. \hfill{$\Box$}
\end{corollary}

This may also be stated as a uniqueness result for $KH$-theory:

\begin{corollary}
Let $k$ be a field of characteristic zero. The homotopy algebraic
$K$-theory $\kh$ is the unique ($\Phi$-rectifiable) functor
$\Sch(k)\longrightarrow \mathrm{Ho}\Spec$, up to equivalence,
which satisfies the descent property (D) and is equivalent to the
algebraic $K$-functor $\K$ over smooth algebraic varieties.
\hfill{$\Box$}
\end{corollary}

\subsection{Algebraic $K$-theory with compact support}
We can apply the same arguments of the proof of theorem
\ref{teorema1} jointly with the compact support extension theorem
in \cite{GN02} to extend the algebraic $K$-theory of smooth
projective varieties over a field of characteristic zero to a
theory with compact support:

\begin{theorem}
Let $k$ be a field of characteristic zero and $\mathbf V(k)$ be
the category of smooth projective $k$-varieties. The rectified
contravariant functor
$$
\K : \mathbf V(k)\longrightarrow {\rm Ho}\Spec
$$
admits a unique extension, up to unique isomorphism of
$\Phi$-rectified functors, to a functor
$$
\K^c : \Sch_c(k)\longrightarrow {\rm Ho} \Spec
$$
such that satisfies the descent property (D) and the compact
support descent property:
\begin{itemize}
\item[($D_c$)] if $Y$ is a subvariety of $X$, then there is a
natural isomorphism
$$
\K^c(X\setminus Y)\cong \holim (\K^c(X)\longrightarrow
\K^c(Y)\longleftarrow \ast).
$$
\end{itemize}
\end{theorem}

In other words, property ($D_c$) says that the sequence
$$
\K^c(X\setminus Y)\lra \K^c(X)\longrightarrow \K^c(Y),
$$
is a fibration sequence in $\rm Ho \Spec$, so that taking homotopy
groups it gives rise to a long exact sequence
$$
\dots\lra K^c_n(X\setminus Y)\lra K^c_n(X)\lra K^c_n(Y)\lra
K^c_{n-1}(X\setminus Y)\lra \dots
$$

In \cite{GS}, theorem $7$, Gillet-Soul{\'e} defined a $\K$-theory with
compact support satisfying ($D_c$), so by the uniqueness of the
compact support extension we find:
\begin{corollary}
Let $X$ be an algebraic variety over a field of characteristic
zero. Then $\K^c(X)$ is naturally isomorphic in $\mathrm{Ho}
\Spec$ to the algebraic $K$-theory with compact support introduced
by Gillet and Soul{\'e} in \cite{GS}, theorem 7. \hfill{$\Box$}
\end{corollary}
We will write $K^c_\ast(X)=\pi_\ast(\K^c(X))$.

\section{Weight filtration}
In this section we prove that there are well defined filtrations
on the groups $KD_\ast (X)$, or equivalently on $KH_\ast (X)$, and
on the groups $K^c_\ast(X)$,  which are trivial for $X$ smooth. In
the compact support case we recover the weight filtration obtained
by Gillet-Soul{\'e}, \cite{GS}.

We fix a field $k$ of characteristic zero.

\subsection{} Let $X$ be an algebraic variety. The spectral sequence
\ref{ssequence} associated to a cubical hyperresolution
$X_\bullet$ of $X$ induces a filtration on the groups $KD_n(X)$.
Our next goal is to prove that this filtration on $KD_{n}(X)$ is
independent of the cubical hyperresolution $X_\bullet$. We will
follow section $3$ of \cite{GN03} closely, where the authors
analyze the weight filtration in an abelian setting.

\subsection{Towers of fibrant spectra} First, we introduce a cohomological
descent structure on the category of towers of fibrations
$\tow(\Spec)$.

\subsubsection{} A tower of fibrations $X(-)$ is a sequence of
fibrations of spectra
$$
\dots \lra X(n)\longrightarrow X(n-1)\lra \dots \lra
X(1)\longrightarrow X(0)\longrightarrow \ast
$$
A morphism of towers of fibrations  is a morphism of diagrams. We
denote by $\tow(\Spec)$ the category of towers of fibrations.

Defining weak equivalences of towers of fibrations and simple
functors for cubical diagrams degree wise, it is immediate to
prove the following result:

\begin{proposition}
The category of towers of fibrations $\tow(\Spec)$ together with
weak equivalences and simple functors for cubical diagrams defined
degree wise is a descent category. {\hfill $\Box$}
\end{proposition}

\subsubsection{} We now introduce a second descent structure on
$\tow(\Spec)$. Recall that if $X(-)$ is a tower of fibrations,
there is a functorial spectral sequence
$$
E_1^{pq}=\pi_{q-p}(F(p)) \Longrightarrow \pi_{q-p}(X),
$$
where $F(p)$ is the fiber of the morphism $X(p)\longrightarrow
X(p-1)$ and $X=\lim X(p)$, see \cite{T80}, 5.43 (where convergence
is understood in the sense of Bousfield-Kan).

\begin{definition}
We say that a morphism of towers $f: X(-)\longrightarrow Y(-)$ is
an $E_2$-\emph{weak equivalence} if the morphism
$E^{\ast\ast}_2(f)$ induced on the $E_2$-terms of the
corresponding spectral sequences is an isomorphism.
\end{definition}

Observe that if $f_p: X(p)\longrightarrow Y(p)$ is a weak
equivalence, for all $p\geq 0$, then $f$ induces an isomorphism in
the $E_1$ terms of the spectral sequence and hence it is also a
$E_2$-weak equivalence.

\subsubsection{} Now we define a simple construction,
$s_2:(\Box , \tow(\Spec))\lra \tow(\Spec)$, compatible
with the $E_2$-weak equivalences: given a tower of fibrations
$X(-)$ and a positive integer $n\geq 0$, we denote by $X[n](-)$
the tower of fibrations defined by
$$
X[n](p):=\left\{
\begin{array}{ll}
  \ast ,& 0\leq p < n, \\
  X(p-n), & p\geq n, \\
\end{array}%
\right.
$$
with the evident morphisms, so that the new tower is obtained by
translating $n$ places to the left the tower $X(-)$.

\begin{definition}\label{diagonal}
Let $\Box$ be a cubical category and $X_\bullet(-)$ be a
$\Box$-diagram of towers of fibrations. Denote by $dX(-)$ the
$\Box$-diagram of towers of fibrations given by
$$
(dX)_\alpha(-) = X_\alpha[|\alpha|-1](-),
$$
with morphisms induced by $X_\bullet$. We define the $s_2$
\emph{simple} of $X_\bullet(-)$ as the tower of fibrations
obtained by applying homotopy limits in each cubical degree of
$dX_\bullet(-)$, that is,
$$
s_2(X_\bullet)(p) := s(dX_\bullet(p)) = \holim_\alpha X_\alpha
(p-|\alpha|+1).
$$
\end{definition}

For example, given a $\Box_1$-diagram $X_\bullet(-)$ of towers of
fibrations
$$
\begin{CD}
\dots @>>> X(1) @>>> X(0) @>>> \ast \\
@. @VVV @VVV @. \\
\dots @>>> Y(1) @>>> Y(0) @>>> \ast \\
@. @AAA @AAA @. \\
\dots @>>> Z(1) @>>> Z(0) @>>> \ast \\
\end{CD}
$$
the new diagram $dX_\bullet(-)$ is the diagram
$$
\begin{CD}
\dots @>>> X(1) @>>> X(0) @>>> \ast \\
@. @VVV @VVV @. \\
\dots @>>> Y(0) @>>> \ast @>>> \ast \\
@. @AAA @AAA @. \\
\dots @>>> Z(1) @>>> Z(0) @>>> \ast \\
\end{CD}
$$
and it follows that its $s_2$ simple in degree $p$ corresponds to
the spectrum
$$
\holim(X(p)\lra Y(p-1)\longleftarrow Z(p)).
$$

\begin{lemma}
For any cubical diagram of towers of fibrations $X_\bullet (-)$
there is a canonical isomorphism of complexes of abelian groups
$$
E^{\ast q}_1(s_2X_\bullet (-)) \lra s(\alpha\mapsto E_1^{\ast
q}(X_\alpha(-))).
$$
\end{lemma}
\begin{proof}
The notation $s(\alpha\mapsto E_1^{\ast q}(X_\alpha(-)))$ refers
to the ordinary simple functor for complexes of abelian groups,
also called the total complex associated to a cubical complex. The
group in degree $p$ of this complex is
$$
s(\alpha\mapsto E_1^{\ast q}(X_\alpha(-)))^p = s(\alpha\mapsto
\pi_{q-r}F_\alpha(r)) =
\bigoplus_{|\alpha|+r=p+1}\pi_{q-r}(F_\alpha (r)),
$$
while the differential is induced by the differentials of the
Bousfield-Kan spectral sequence of the tower $X_\alpha(-)$.

On the other hand, by definition, for each $p$, $s_2X_\bullet (p)$
is the ordinary simple of the cubical diagram of spectra
$dX_\bullet(p)$, so the complex $E^{\ast q}_1(s_2X_\bullet (-))$
is the $E_1$-term of the Bousfield-Kan spectral sequence
associated to the tower of fibrations
$$
\dots\lra sdX_\bullet (p)\lra \dots\lra sdX_\bullet (1)\lra
sdX_\bullet(0)\lra \ast
$$
Denote by $F_\alpha(p)$ the fiber of the fibration $ X_\alpha
(p)\lra X_\alpha (p-1)$. Since homotopy limits commute, the fiber
of the fibration $sdX_\bullet (p)\lra sdX_\bullet (p-1)$ is
isomorphic to the simple spectrum associated to the cubical
diagram $dF_\bullet(p)$. But, in this diagram all morphisms are
constant, so
$$
sdF_\bullet(p) = \prod_{\alpha}\Omega^{|\alpha |
-1}F_\alpha(p-|\alpha |+1) =
\prod_{|\alpha|+r=p+1}\Omega^{p-r}F_\alpha(r),
$$
hence its homotopy groups are given by
$$
E_1^{pq} = \pi_{q-p}(sdF_\bullet(p)) =
\bigoplus_{|\alpha|+r=p+1}\pi_{q-r}(F_\alpha (r)).
$$
The differential is also induced by the differentials of the
Bousfield-Kan spectral sequence of the tower $X_\alpha(-)$.
\end{proof}

\begin{proposition}
The simple $s_2$ and the $E_2$-weak equivalences define a
cohomological descent category structure on $\tow(\Spec)$.
\end{proposition}
\begin{proof}
Observe that a morphism between towers of fibrations $f$ is a
$E_2$-weak equivalence if and only if the morphism $E_1(f)$ of the
corresponding spectral sequence is a quasi-isomorphism of
complexes. If $GrC_\ast(\mathbb Z)$ denotes the category of graded
complexes of abelian groups, the functor
\begin{eqnarray*}
E_1: \tow(\Spec)&\lra & GrC_\ast(\mathbb Z)\\
X(-) &\longmapsto & E_1^\ast
\end{eqnarray*}
commutes with direct sums and, by the previous result, it commutes
with the simple $s_2$ functor, so the result follows from
\cite{GN02}, $(1.5.12)$.
\end{proof}

\subsection{An extension criterion for towers}
In the next result we write $\mathrm{Ho}_2(\tow(\Spec))$ for the
homotopy category obtained from $\tow(\Spec)$ by inverting
$E_2$-weak equivalences.

The following result, remarked by Navarro several years ago in the
abelian context, is the key point in order to extend some functors
on $\Reg(k)$ with values in the category of spectra to functors
defined for all varieties and taking values in
$\mathrm{Ho}_2(\tow(\Spec))$.

\begin{proposition}\label{descensfilt}
{\em [Compare with \cite{GN03}, proposition (3.10)].} Let
$G:\Reg(k)\longrightarrow \mathrm{Ho}\Spec$ be a
$\Phi$-rectifiable functor and denote also by
$$G:\Reg(k)\longrightarrow \mathrm{Ho}_2(\tow(\Spec))$$
the associated constant functor. Then, $G$ satisfies property
$(F2)$ if and only if for every elementary acyclic square the
sequence
$$
0\longrightarrow \pi_nG(X)
\stackrel{f^\ast-i^\ast}{\longrightarrow} \pi_nG(\widetilde
X)\oplus \pi_nG(Y)\stackrel{j^\ast + g^\ast}{\longrightarrow}
\pi_nG(\widetilde Y)\longrightarrow 0,
$$
is exact.
\end{proposition}
\begin{proof}
The $(F2)$ property for the extended functor $G$ says that the
morphism
$$
E_1^{\ast ,q}(G(X))\lra E_1^{\ast ,q}s_2 G(X_\bullet)
$$
is a quasi-isomorphism. Observe that we have
$$
E_1^{\ast ,q}(G(X),\tau) = \left\{
\begin{array}{ll}
  \pi_q(G(X)), & p=0, \\
  0, & p\not= 0. \\
\end{array}%
\right.
$$
By the other hand, the $E_1$-page of the spectral sequence of $s_2
G(X_\bullet)$ reduces to the exact sequence
$$
\pi_nG(\widetilde X)\oplus \pi_nG(Y)\stackrel{j^\ast +
g^\ast}{\longrightarrow} \pi_nG(\widetilde Y),
$$
so the $(F2)$ property is equivalent to the the fact that the
morphism of complexes of abelian groups
$$
\pi_nG(X)
\stackrel{f^\ast-i^\ast}{\longrightarrow}(\pi_nG(\widetilde
X)\oplus \pi_nG(Y)\stackrel{j^\ast + g^\ast}{\longrightarrow}
\pi_nG(\widetilde Y)),
$$
is a quasi-isomorphism, which is precisely the condition stated in
the proposition.
\end{proof}

\subsubsection{} We return now to the applications to algebraic
$K$-theory. The following proposition has also been proved by
Gillet-Soul{\'e} directly from Thomason's calculations, see \cite{GS},
theorem 5:

\begin{proposition}\label{exacta}
For any elementary acyclic square of $\Reg(k)$ and any $n\geq 0$,
the sequence
$$
0\longrightarrow K_n(X) \stackrel{f^\ast-i^\ast}{\longrightarrow}
K_n(\widetilde X)\oplus K_n(Y)\stackrel{j^\ast +
g^\ast}{\longrightarrow} K_n(\widetilde Y)\longrightarrow 0,
$$
is exact.
\end{proposition}
\begin{proof}
As we have recalled in the proof of theorem \ref{teorema1}, an
elementary acyclic square gives rise to a homotopy cartesian
square of algebraic $K$-theory spectra, so we have an exact
sequence
$$
\dots \lra K_n(X)\stackrel{f^\ast -i^\ast}{\lra }K_n(\widetilde
X)\oplus K_n( Y)\stackrel{j^\ast + g^\ast}{\lra }K_n({\widetilde
Y})\stackrel{\delta}{\lra }K_{n-1}(X)\lra \dots
$$
But, by Thomason calculation of the algebraic $K$-theory of a blow
up (\cite{T93}), there are isomorphisms
\begin{eqnarray*}
\varphi : K_n(X)\bigoplus_{i=1}^{d-1}K_n(Y) &\longrightarrow &
K_n(\widetilde X),\\
\psi : \bigoplus_{i=0}^{d-1}K_n(Y) &\longrightarrow &
K_n(\widetilde Y),
\end{eqnarray*}
given, respectively,  by
\begin{eqnarray*}
\varphi (x,y_1,\dots ,y_{d-1}) &=& f^\ast (x) +
\bigoplus_{i=1}^{d-1}
j_\ast(\ell^{-i}\cup g^\ast(y_i)),\\
\psi (y_0,y_1,\dots ,y_{d-1}) &=& \bigoplus_{i=0}^{d-1}
j_\ast(\ell^{-i}\cup g^\ast(y_i)).
\end{eqnarray*}
With this identifications the morphism $f^\ast$ corresponds to the
inclusion of $K_n(X)$ on the first factor of $K_n(\widetilde X)$,
and so the morphism $f^\ast -i^\ast$ is injective. This splits the
exact sequence above into the required short exact sequences.
\end{proof}

Now, by propositions \ref{descensfilt} and \ref{exacta} we  can
apply the extension criterion of theorem \ref{extensio}, so we
find:

\begin{corollary}
Let $k$ be a field of characteristic zero. The constant algebraic
$K$-theory functor $\K :\Reg(k)\longrightarrow
\mathrm{Ho}_2(\tow(\Spec))$ admits an essentially unique extension
$\kd(-):\Sch(k)\longrightarrow \mathrm{Ho}_2(\tow(\Spec))$ which
satisfies the descent property (D). Moreover, for any variety $X$,
the tower of fibrations $\kd(-)(X)$ satisfies $\kd(n)(X) = \kd(X)$
for $n\gg 0$.
\end{corollary}
\begin{proof}
We have only to justify the last sentence. Take an algebraic
variety $X$ and an hyperresolution $X_\bullet$, whose type $\Box$
is of length $\ell$. By the definition of the descent functor
$\kd(-)$, the tower $\kd(-)(X)$ is the $s_2$-simple tower
associated to the diagram of constant towers $\K(X_\bullet )$,
that is, it is the tower whose spectra are the homotopy limits of
the diagram $d\K(X_\bullet)(n)$ for each $n$ (see definition
\ref{diagonal}). Observe that this diagram is constant for $n\geq
\ell$ and, moreover, it is precisely the cubical diagram
$X_\bullet$, so the result follows.
\end{proof}

Since the spectral sequence of a tower of fibrations is functorial
in the category $\mathrm{Ho}_2(\tow(\Spec))$ from the $E_2$-term
on, we deduce from the corollary above:

\begin{corollary}
There is a well defined and functorial finite increasing
filtration $F^p$ on $KD_n(X)$ which is trivial for smooth
varieties.
\end{corollary}

\begin{remark}
Equivalently, by \ref{kh=kd}, for any variety $X$ the last
corollary defines a functorial finite filtration on the homotopy
algebraic $K$-theory groups $KH_n(X)$.
\end{remark}

\subsection{} Finally, we observe that the same procedure may be applied
to the algebraic $K$-theory with compact support. In this case,
from the uniqueness property of descent extensions and \cite{GS},
theorem 7, we deduce:

\begin{corollary}
There is a well defined and functorial finite increasing
filtration $W^pK^c_n(X)$ which is trivial for complete smooth
varieties. This filtration coincides with the weight filtration
defined by Gillet-Soul{\'e} in \cite{GS}.
\end{corollary}

\end{large}
\end{document}